\documentclass[11pt]{amsart}%
\usepackage{amstext}
\usepackage{amsthm,euscript}
\usepackage{amscd}
\usepackage{amssymb}
\usepackage{amsmath}
\usepackage{mathrsfs}% for mathscr %%{MnSymbol} %for \ostar
\usepackage{mathtools}
\usepackage[all]{xy}
\usepackage{paralist}%% list packages
\usepackage{enumitem}%% for change of indentation of enumeration,
\usepackage{moreenum}
\usepackage{scalerel,stackengine}%%for reallywidehat
\usepackage[colorlinks=true, pdfstartview=FitV, linkcolor=blue,
citecolor=blue,urlcolor=blue,breaklinks=true]{hyperref}
\usepackage{subfiles}

\swapnumbers 

\newcommand{\lv}[1]{} %% use this if the argument of \lv should not appear
\newcommand{\comments}[1]{}
\newcommand{\pcomments}[1]{}
\newcommand{\inparcom}[1]{}
\newcounter{question} 
\newcommand{\quest}[1]{} %%% if questions will not appear in the pdf
\newcommand{\arxiv}[1]{\href{http://arxiv.org/abs/#1}{\tt arXiv:\nolinkurl{#1}}}
\newtheorem{prop}[subsection]{Proposition}

\newtheorem{lem}[subsection]{Lemma}

\newtheorem{cor}[subsection]{Corollary}

\theoremstyle{definition}

\newtheorem{remarks}[subsection]{Remarks}

\newtheorem{leer}[subsection]{}%

\numberwithin{equation}{subsection}%% was {section} before

\newcommand\sm{\smallskip}
\newcommand\ms{\medskip}

\def\limind{\mathop{\oalign{lim\cr\hidewidth$\longrightarrow$\hidewidth \cr}}}%directed limit
\def\limproj{\mathop{\oalign{lim\cr\hidewidth$\longleftarrow$\hidewidth\cr}}}

\newcommand{\ol}{\overline} 

\newcommand{\simlgr}{\buildrel \sim \over \longrightarrow}

\newcommand{\wdh}{\widehat}
\newcommand{\longto}{\longrightarrow}

\def\co{\colon}
\def\ot{\otimes}

\newcommand{\me}{^{-1}}
\def\dar[#1]{\ar@<2pt>[#1]\ar@<-2pt>[#1]}
\def\wtl{\widetilde}

\stackMath
\newcommand\reallywidehat[1]{%
\savestack{\tmpbox}{\stretchto{%
  \scaleto{%
    \scalerel*[\widthof{\ensuremath{#1}}]{\kern.1pt\mathchar"0362\kern.1pt}%
    {\rule{0ex}{\textheight}}%WIDTH-LIMITED CIRCUMFLEX
  }{\textheight}%
}{2.4ex}}%
\stackon[-6.9pt]{#1}{\tmpbox}%
}

\newcommand{\can}{\operatorname{can}}
\newcommand{\Frac}{\operatorname{Frac}}

\newcommand{\Id}{\operatorname{Id}}
\newcommand{\Ima}{\operatorname{Im}}

\newcommand{\inc}{{\operatorname{inc}}}

\newcommand{\Ker}{\operatorname{Ker}}

\newcommand{\Spec}{\operatorname{Spec}}
\newcommand{\calO}{\operatorname{\mathcal O}}

\newcommand{\scJ}{\mathscr{J}}
\newcommand{\scK}{\mathscr{K}}

\newcommand{\scO}{\mathscr{O}}

\newcommand{\m}{\mathfrak m} %maximal ideal

\newcommand{\kalg}{k\mathchar45\mathbf{alg}}

\newcommand{\Ralg}{R\mathchar45\mathbf{alg}}

\newcommand\al{\alpha}

\newcommand\ka{\kappa}
\newcommand\la{\lambda} 
 \newcommand\vphi{\varphi}

\title{Surjectivity and flatness over DVR's 
(after Moret-Bailly)}

\author[B.~Margaux]{Benedictus Margaux}
\address{Laboratoire de Recherche ``Princesse St\'ephanie'', Monte Carlo 51840, Monaco}
\email{benedictus.margaux@gmail.com}
\date{\today}

\begin{document}
\begin{abstract} We study morphisms of schemes $f \co X \to S$ which are locally of finite type. We present conditions under which there exists a
morphism $g \co S'\to X$ of $S$--schemes such that $f \circ g $ is the canonical morphism $S'\to S$. Furthermore, we exhibit situations in which $f$ is flat surjective. Our results are mostly concerned with $S$ being the spectrum of a DVR. 
\end{abstract}

\maketitle

\bigskip


\smallskip

\bigskip
\bigskip

\section*{Introduction}
This note presents criteria which establish surjectivity and flatness of a scheme morphism $f\co X \to S$ which is locally of finite type. 

Section~\S\ref{sec:weaksurjectivity} is devoted to a weak version of surjectivity. Given a noetherian local integral domain $R$ and a faithfully flat (= flat and surjective) morphism $f\co X \to \Spec(R)$ which is locally of finite type, we show in \ref{lem_wsur} that there exists a valuation ring $B$ dominating $R$ and a morphism $g\co \Spec(B) \to X$ such that $f\circ g$ is the canonical map $\Spec(B) \to \Spec(R)$, i.e., the diagram 
\[ \xymatrix{
  \Spec(B) \ar[rr]^g \ar[dr] && X \ar[dl]^f \\ & \Spec(R)
}\]
commutes. Imposing more conditions on $R$ allows one to make the structure of $B$ more precise, e.g., if $R$ is a DVR we can assume that $B$ is a DVR. As the picture indicates, one can view the morphism $g$ as a weak form of a section of $f$. Indeed, our proof of \ref{lem_wsur} is based on a result from \cite[IV$_3$]{EGA} regarding morphisms with everywhere local sections, which we review in \ref{elos}.%% 

We consider a similar situation in section~\S\ref{sec:ss}: $S$ is a locally noetherian integral scheme and $f \co X \to S$ is a morphism locally of finite type. We are interested in conditions ensuring that $f$ is flat surjective. In this generality, we assume that there exists a faithfully flat 
  %surjective 
morphism $S'\to S$ and an $S$--morphism $S'\to X$ which factors $S'\to S$. Denoting by $X_\eta$ the generic fibre of $f$ and by $\wtl X = \Ima (X_\eta \to X)$ the schematic image (= schematic closure) of $X_\eta$, our vehicle in \S\ref{sec:ss} is the induced morphism $\wtl f \co \wtl X \to S$:  
\[ \xymatrix@C=50pt{
  S' \ar[r] \ar[dr] & X \ar[d]^f & \wtl X \ar[l]\ar[dl]^{\wtl f} \\ & S
}\]
which we first assume and later show in special cases to be flat. 
The problem then becomes lifting flatness from $\wtl f$ to $f$. We address this in decreasing generality in \ref{MBco}--\ref{prop_strong_surj}, ending with the case of $S$ being the spectrum of a henselian DVR. Crucial for this is a characterization of the points in $\wtl X$, presented in \ref{MB}\eqref{MBa}. 
Finally, in Proposition~\ref{AGB1}, we investigate the situation where $S$ is a Dedekind scheme, not necessarily irreducible. 
\ms

{\em Acknowledgements.} The results in this note are essentially all due to Laurent Moret-Bailly. We thank him for encouraging us to write them up. We also thank   
Qing Liu for helpful comments on an earlier version of this note and Marion Jeannin for her interest in this paper.

\section{Weak Surjectivity for morphisms}\label{sec:weaksurjectivity}

\begin{leer}{\bf Morphisms with everywhere local sections.}
\label{elos} A morphism $f\co X \to S$ has {\em everywhere local sections\/} (with respect to the Zariski topology), if for every $s\in S$ there exists an open neighbourhood $U$ of $s$ and a morphism $s_U \co U \to X$ such that $f \circ s_U = \Id_U$,  equivalently, there exists a Zariski cover $\{U_i\}$ of $S$ and sections of $f$ over every $U_i$. 
We will use this concept in the proof of the following Lemma~\ref{lem_wsur}. But first  we mention some useful facts:
\sm

\begin{inparaenum}[(a)] \item \label{elos-aa} ({\em $S$ local}) If $S=\Spec(R)$ for a local ring $(R,\m)$ and $f$ has everywhere local sections, then $f$ has a (global) section given that $S$ is the only open neighbourhood of the closed point $\m\in S$. \sm

\item \label{elos-a} ({\em Base change}) If $g \co S' \to S$ is a morphism of schemes and if $f \co X \to S$ has everywhere local sections, then the same holds for the morphism $f' \co X'=X\times_S S' \to S'$.

%\inparcom{(2025-01-27) Put justification for base change in the long version}
\lv{%%%%%%%%%%%%%%%   lv starts here
Indeed, fix $s'\in S'$ and let $U\subset S$ be an open neighbourhood of $g(s)$ for which there exists a section $s_U \co U \to X$. Then $U' = g\me(U)$ is an open neighbourhood of $s'$, and we have the inclusion morphism $\inc_{U'} \co U' \to S'$ as well as the morphism $s_U \circ g \circ \inc_{U'} \co U' \to X$ satisfying
$f \circ (s_U \circ g \circ \inc_{U'}) = g \circ \inc_{U'}$.
\[ \xymatrix@C=50pt{
 X'  \ar[rrr]^{g'} \ar[dd]_{f'}&&&  X \ar[dd]^f
 \\
 &  U' \ar[r]^{g \circ \inc_{U'}} \ar@{-->}[ul]_{\exists!\; s_{U'}}  \ar[dl]_{\inc_{U'}} & U \ar[ur]^{s_U} \ar[dr]^{\inc_U}
 \\
S' \ar[rrr]^g &&& S
 }\]
Hence there exists a unique morphism $s_{U'} \co U' \to X$ such that $f' \circ s_{U'} = \inc_{U'}$, i.e., $s_{U'}$ is a local section of $f'$.} %%%%%%%%%%%   end of lv 

\item \label{elos-b} (\cite[IV$_3$, 14.5.10]{EGA}) Let $S$ be a noetherian scheme, and let $f \co X \to S$ be a morphism that is universally open, locally of finite type and surjective. Then there exists a finite surjective morphism $g \co S' \to S$ such that the morphism $f' \co X\times_S S' \to S'$ has everywhere local sections.
 \end{inparaenum}
\end{leer}
 
%\comments{(2025-01-27) Deleted the following remark  since \cite[Tag 02F2]{Stacks} is %wrong, as Laurent claims: 
%{\tt % \item \label{elos-d} 
%A related situation is discussed in \cite[Tag 02F2]{Stacks}: if $X\to S$ is a morphism all of whose fibers are nonempty, does there exists a finite surjective morphism $S'\to S$ such that the base change $X_{S'} \to S'$ has a global section? This is true if $S=\Spec(A)$, $A$ a Dedekind ring with finite residue fields at closed points, and the morphism $X\to S$ is flat with geometrically irreducible fibre.}} 

\begin{prop}[Weak surjectivity] \label{lem_wsur} Let $R$ be a noetherian local integral domain  and let $f \co X \to \Spec(R)$ be a faithfully flat 
%% recall faithfully flat = flat and surjective 
morphism which is locally of finite type. Then there exists a valuation ring $B\in \Ralg$ with the following properties:
\begin{enumerate}[label={\rm (\roman*)}]
  \item \label{lem_wsur-i} there exists a morphism $g \co \Spec(B) \to X$ such that $f\circ g$ is the canonical morphism $\Spec(B) \to \Spec(R)$;

  \item \label{lem_wsur-i-i} $B$ dominates $R$ in the sense that $R\subset B$ and $\m_B \cap R = \m_R$, where $\m_B$ and $\m_R$ are the maximal ideals of $B$ and $R$ respectively;

  \item \label{lem_wsur-iii} the fraction field $\Frac(B)$ of $B$ is a finite extension of $\Frac(R)=K$.

\end{enumerate}
Moreover, in the following special cases for $R$ listed below, the algebra $B$ can be assumed to have additional properties besides \ref{lem_wsur-i}--\ref{lem_wsur-iii}:\sm

\begin{inparaenum}[\rm (a)] \noindent \item \label{lem_wsur-a}
If $R$ is a DVR, we can assume that \end{inparaenum}
\begin{enumerate}[label={\rm (\roman*)}] \setcounter{enumi}{3}

\item \label{lem_wsur-iv} $B$ is a DVR and a %\label{lem_wsur-iv} $B$ is a
   faithfully flat $R$--module.
\end{enumerate}
\sm

\begin{inparaenum}[\rm (a)] \setcounter{enumi}{1}
\noindent \item \label{lem_wsur-b} If $R$ is a henselian DVR, then besides \ref{lem_wsur-i}--\ref{lem_wsur-iv} we can suppose that \end{inparaenum}
\begin{enumerate}[label={\rm (\roman*)}] \setcounter{enumi}{4}

\item \label{lem_wsur-iv'} $B$ is a henselian DVR and the integral closure of $R$ in a finite field extension of $K$. 
    % hence in particular a normal ring.
\end{enumerate}
\sm

\begin{inparaenum}[\rm (a)] \setcounter{enumi}{2}
\noindent \item \label{lem_wsur-c}
Finally, if $R$ is a henselian DVR and a japanese ring, then \ref{lem_wsur-i}--\ref{lem_wsur-iv'} hold as well as \end{inparaenum}
\begin{enumerate}[label={\rm (\roman*)}]\setcounter{enumi}{5}
\item\label{lem_wsur-v} $B$ is finite over $R$.
\end{enumerate}
\end{prop}

\begin{proof} Since $S=\Spec(R)$ is a noetherian scheme, %\cite[Prop. 3.19]{GW}
$f$ is locally of finite presentation %\cite[10.36]{GW}
and therefore universally open by \cite[Tag 01UA]{Stacks}. %\cite[14.33]{GW}.
We can thus apply \ref{elos}\eqref{elos-b}. %\cite[IV$_3$, 14.5.10]{EGA}
The morphism $g\co S' \to S$, whose existence is established there, is finite, hence %affine, so that 
$S' = \Spec(R')$ for a finite $R$--algebra $R'$. Moreover, the morphism $f'\co X'=X \times_S S' \to S'$ has everywhere local sections.
%\comments{(2022-01-24) Until here we need: $S$ noetherian, i.e., $R$ noetherian, and $f$ locally of finite type, surjective + universally open. We do not need faithfully flat, i.e., flat + surjective. But we use faithfully flat, in order to apply \cite[01UA]{Stacks}. }
%\pcomments{(2022-0127) I agree with you that we could replace flat by universally open which is very close. Actually, if the fiber of $f$ at the closed point is geometrically reduced, this is the same thing according to \cite[IV$_3$, Cor.~15.2.3]{EGA}.}

 Let $K=\Frac(R)$ be the fraction field of the integral domain $R$, and let $L$ be  the  total ring of fractions of the finite-dimensional $K$--algebra $R'\ot_R K$ \cite[Tag 02C5]{Stacks}. 
 %\footnote{AP. What is ``a quotient field" of a $K-$algebra?}
%\lv{ this exists since any non-zero element of K in R'\ot_R K is invertible, so that any maximal ideal of $R'\ot_R K$ intersects $K$ trivially} 
Then $L$ is a finite extension of $K$ and $R$ is a subring of $L$. By \cite[VI, \S1.3, Thm.~3]{BAC2} the integral closure $A$ of $R$ in $L$ is the intersection of
all valuation rings of $L$ dominating $R$,  in particular we can choose a valuation ring $B$ of $L$  that dominates $R$, i.e., $A\subset B$. Note then that \ref{lem_wsur-iii} also holds because $\Frac(B) = L$. Since $R'$ is an integral extension of $R$, the image of the canonical homomorphism $R'\to R'\ot_R K \to L$ is an integral extension of $R$, thus contained in $A$, hence in $B$, and in this way induces a homomorphism $R'\to B$. Thus we have the following situation ($\inc = $ inclusion):
\[ \xymatrix@C=50pt{
  R \ar[d]_{\inc}  \ar[r]& R' \ar[r] \ar[dl] \ar@{-->}[d]& R'\ot_R K \ar@{>>}[d]
  \\
  A  \ar[r]^{\inc}  & B\;  \ar[r]^{\inc} & L
  }\]
The homomorphism $R'\to B$ gives rise to a morphism $\Spec(B) \to S'=\Spec(R')$, which we use as base change map for the morphism $f'$, thus giving rise to a morphism
\[
     f''\co X\ot_{S} \Spec(B) = X'\ot_{S'} \Spec(B) \longto \Spec(B)
\] of schemes.
By \ref{elos}\eqref{elos-a}  we know that $f''$ admits everywhere local sections. But since $B$ is a local ring, $f''$ has in fact a global section by  \ref{elos}\eqref{elos-aa}. Composing this section with the canonical morphisms $X'\times_{S'} \Spec(B) \to X'= X\times_S S' \to X$ shows $X(B) \ne \emptyset$, i.e., \ref{lem_wsur-i} is true.
\sm

\eqref{lem_wsur-a} Suppose now that $R$ is a DVR, thus a noetherian valuation ring \cite[VI, \S3.6, Prop.~9]{BAC2} and, in particular, a local integral domain.
We can therefore apply the above
% and \cite[VI, \S1.3, Cor.~3]{BAC2}, to conclude that $A$ is the intersection of all valuation rings of $L$ dominating $R$.
and \cite[VI, \S8.1, Cor.~3]{BAC2} to conclude that $B$ is a discrete valuation ring. Also, $B$ is a faithful and torsion-free $R$--module, proving \ref{lem_wsur-iv} by applying \cite[VI, \S3.6, Lem.~1]{BAC2}.
\sm

\eqref{lem_wsur-b} If $R$ is a henselian DVR, the discrete valuation of $K$ extends uniquely to $L$ \cite[32.8]{Warner}, implying $A=B$. It is a henselian DVR by \cite[32.12]{Warner}. %icular, \ref{lem_wsur-iv'} holds.
%% Other reference for unique extension is Neukirch, Algebraic Number Theory, Ch II, 6.2
\sm

\eqref{lem_wsur-c} Since $A=B$, the claim holds by the definition of a japanese ring
\cite[0; 23.1]{EGA} (or \cite[IX, \S4.1]{BAC3}, \cite[Tag 032F]{Stacks}): a noetherian integral domain $R$ is a japanese ring if the integral closure of $R$ in every finite extension field of $\Frac(R)$ is a finite $R$--algebra. \end{proof}
\sm

Instead of assuming  in \ref{lem_wsur}\eqref{lem_wsur-c} that $R$ is a japanese henselian DVR, we could have assumed that $R$ is an excellent henselian DVR, because of the following result:

\begin{prop}\label{Sch} Let $R$ be a DVR with fraction field $K$. We denote by $\wdh K$ the completion of $K$ with respect to the discrete valuation of $K$. Then the following are equivalent:

\begin{enumerate}[label={\rm (\roman*)}]

         \item \label{excellent_i} $R$ is excellent;

         \item \label{excellent_ii} $\wdh K$ is separable over $K$;

         \item \label{excellent_iii} $R$ is japanese.

       \end{enumerate}
\end{prop}

\begin{proof}
  For \ref{excellent_i} $\iff$ \ref{excellent_iii} see \cite[IV$_2$, 7.8.3]{EGA}, and for \ref{excellent_ii} $\iff$ \ref{excellent_iii} see \cite[IX, \S4.4, Thm.~3]{BAC3}. The result is also proven in \cite[Prop~4.1]{Sc20}. \end{proof}

%\pcomments{Remark (a) after G3. I think that the point of the proposition is that a japanese DVR is excellent.}

\begin{remarks} \label{excex}\begin{inparaenum}[(a)]
  \item \label{Sc20-a} The point of Proposition~\ref{Sch}  is that a japanese DVR is   excellent; the converse holds in general: any excellent ring is japanese \cite[IV$_2$, 7.8.3(vi)]{EGA}.

Let us give a quick proof of \ref{excellent_ii} $\implies$ \ref{excellent_iii}. We denote by $\wdh R$ the closure of $R$ in $\wdh K$. Recall that $\wdh R =  \limproj R/\pi^n R$, where $\pi$ is an uniformizing parameter \cite[\S II.1]{Se62}. One also knows that the canonical map $ K \otimes_R \widehat R \simlgr \wdh K$
is an isomorphism. The assumption \ref{excellent_ii} implies that the $K$--algebra
$K \otimes_R \widehat R$ is separable, so that $R$ is japanese by \cite[IX, \S4.2, Prop.~3]{BAC3}.
\sm

\item \label{excex-b} Let $R$ be a DVR with fraction field $K$. If $K$ has characteristic $0$, then $R$ is excellent (\cite[Prop.~3.1]{RL}, \cite[Tag 07QW]{Stacks}). This fails if $K$ is of characteristic $p > 0$, see \cite[\S 11.5]{RL}. On the other hand, examples of excellent DVR's are obtained from Kunz's result \cite[Th. 2.5]{Ku}: If $K$ is of characteristic $p>0$ and $R$ is finite over $R^p$, then $R$ is excellent. 
\end{inparaenum}
\end{remarks}

\section{Strong surjectivity in the DVR case}\label{sec:ss}

\begin{leer}{\bf Schematic image.}\label{sim} Let $f\co X \to Y$ be a morphism of schemes. There exists a unique closed subscheme $\Ima(f)$ of $Y$, called the {\em schematic image},  satisfying the following two conditions:\begin{enumerate}[label={\rm (\roman*)}]
  \item\label{simi} $f$ factors through the inclusion $\inc_{\Ima(f)} \co \Ima(f)  \hookrightarrow  Y$, and 
  
  \item $\Ima(f)$ is the ``smallest'' closed subscheme with property \ref{simi}: if $Z\subset Y$ is a closed subscheme of $Y$ and $f$ factors through the inclusion $\inc_Z \co Z \hookrightarrow Y$, then $\inc_{\Ima(f)}$ factors through $\inc_Z$: 
  \[ \xymatrix@C=70pt{ X \ar[r]^f \ar[d]   &Y 
     \\ \Ima(f) \ar[r] \ar[ur]_{\inc_{\Ima(f)}} & Z\ar[u]_{\inc_Z}     
  }
  \]     
\end{enumerate} 
(\cite[Lem.~10.29]{GW}, or \cite[Tag 01R7]{Stacks} where $\Ima(f)$ is called the {\em scheme theoretic image of $f$}). In general, the schematic image has some strange properties, see \cite[Tag 0GIK]{Stacks}, but is well-behaved for quasi-compact maps:%
\sm 

\begin{inparaenum}[(a)] 
\item \label{sima} (\cite[Prop.~10.30]{GW}, \cite[Tag 01R8]{Stacks}) Let $f \co X \to Y$ be a quasi-compact morphism. %cite[(10.1)]{GW}
 Then $\scK_f = \Ker(f^\flat \co \scO_Y \to f_* \scO_X)$ is a quasi-coherent ideal of $\scO_Y$, $\Ima(f) = \Spec(\scO_Y/\scK_f)$ is the closed subscheme of $Y$ determined by $\scK_f$, and the underlying topological space of $\Ima(f)$ is the Zariski closure of $f(X)$ in $Y$.  \sm 
 
\item \label{simb} ({\em Flat base change} \cite[Lem.~14.6]{GW}) Let $f \co X \to Y$ be a quasi-compact morphism and let $S'\to S$ be a {\em flat\/}  morphism. Then $\Ima(f) \times_S S' = \Ima(f')$ for the base change $f'\co X \times_S S' \to Y\times_S S'$. \sm

\item \label{simc} ({\em Functoriality} \cite[Tag 01R9]{Stacks}) Let 
\begin{equation}\label{simc1} 
\vcenter{ \xymatrix@C=70pt{X_1 \ar[r]^{f_1} \ar[d] &X \ar[d]^h 
   \\ Y_1 \ar[r]^{g_1}  & Y 
}}\end{equation}          
be a commutative diagram of schemes. Then there exists a unique morphism $h_1\co \Ima(f_1) \to \Ima(g_1)$ such that the two squares in the diagram below commute: 
\begin{equation}\label{simc2} \vcenter{
\xymatrix@C=45pt{X_1 \ar[r]^{f_1'} \ar[d] &   \Ima(f_1)\ar[d]^{h_1} \ar[r]^\inc  &X \ar[d]^h 
   \\ Y_1 \ar[r]^{g_1'} &\Ima(g_1) \ar[r]^\inc   & Y 
}}\end{equation}
where $f_1'$ and $g_1'$ are the maps given by $f_1$ and $g_1$ respectively. \sm

\item \label{simex} {\em Example: Schematic closure of the generic fibre.} Let $S$ be an irreducible scheme with generic point $\eta$, and let $h\co Z \to X$ be a morphism of schemes over $S$. Denoting by $Z_\eta$ and $X_\eta$ the generic fibres of $Z$ and $X$ respectively, we get a commutative diagram as in \eqref{simc1}, hence also the commutative diagram \eqref{simc2}: 
\begin{equation}\label{simex1} 
\vcenter{ \xymatrix@C=40pt{Z_\eta \ar[r]^{j_X} \ar[d]_{h_\eta} 
  &   Z \ar[d]^h 
  \\ X_\eta \ar[r]^{j_X}  & X }}
\quad \rightsquigarrow \quad
\vcenter{
\xymatrix@C=30pt{Z_\eta \ar[r]^{j_Z} \ar[d]_{h_\eta} 
  &\Ima(j_Z)\ar[d]^{\wtl h} \ar[r]^{i_Z}  &
   Z \ar[d]^h 
 \\ X_\eta \ar[r]^{j_X} & \Ima(j_X) \ar[r]^{i_X} & X 
}}\end{equation}
In this situation the schematic image $\Ima(j_Z)$ is also referred to as the {\em schematic closure of $Z_\eta$}. 
\sm 

In this setting let $X_S(Z)$ denote the set of $S$--morphisms $Z \to X$, put $\wtl X = \Ima(j_X)$, and assume that $Z=\Ima(j_Z)$. Then  
\begin{equation}\label{simex2}
  \wtl X_S(Z) \simlgr X_S(Z), \quad \al \mapsto i_X \circ \al
  \end{equation} 
is a bijection.  Indeed, the map is well-defined and injective. Given $h\in X_S(Z)$, we apply \eqref{simex1} and get $\wtl h$ satisfying $ h = h \circ i_Z = i_X \circ \wtl h$.

In \eqref{simEGA} we will give an example where the assumption $Z=\Ima(j_Z)$ is fulfilled.\sm

\item\label{simEGA} {\em Let $S$ be a locally noetherian regular irreducible scheme of dimension $\le 1$ and let $f \co X \to S$ be a morphism. Then the induced morphism $\wtl f \co \Ima(j) \to S$ is flat, where $\Ima(j)$ is the schematic closure of the generic fibre of $X$. Moreover, 
\begin{equation}
  \label{simEGA1}  \text{$f$ is flat} \quad \iff \quad \Ima(j) = X. 
\end{equation} } 
Indeed,  if $\dim S = 0$, then $S=\Spec(k)$ for $k$ a regular local ring of dimension $0$. Thus, $k$ is a field and flatness is clear. %[GW; Remark B.78(2)] 
If $\dim S = 1$, then flatness of $\wtl f$ follows from  \cite[IV$_2$, (2.8.3)]{EGA}, using that by {\em loc.~cit.} for $S=Y$ and $Z'=X_\eta$, there exists a unique closed subscheme $X^+$ which is flat over $S$; it is given as $X^+ = \Spec(\scO_X/ \scJ)$ where $\scJ = \Ker(j^\flat \co \scO_X \to \scO_{X_\eta})$ and therefore equals $\Ima(j)$ by \ref{sim}\eqref{sima}. 

If $f$ is flat, then $\Ima(j) = X$ follows from the uniqueness assertion in \cite[IV$_2$, (2.8.5)]{EGA}. Conversely, if $\Ima(j)=X$, then $\wtl f = f$ is flat.  
\end{inparaenum}
\end{leer}

\sm 

The following Lemma~\ref{MB} is a preparation to Proposition~\ref{MBco}, where we will combine its parts.

\begin{lem}\label{MB} Let $S$ be an integral scheme, and let $f\co X \to S$ be a morphism that is locally of finite type. We denote the fibre of $f$ over $s\in S$ by $X_s$, the generic point of $S$ by $\eta$,  the generic fibre of $f$ by $X_\eta$, the schematic image of the canonical morphism $j \co X_\eta \to X$ by $\wtl X$, let $i \co \wtl X \to X$ be the canonical morphism and for $s\in S$ let $i_s \co \wtl X_s \to X_s$ be the base change of $i$ which, we recall, is a closed immersion as is $i$:  
\begin{equation}\label{MB0} \vcenter{
\xymatrix@C=45pt{
  & X \ar[d]^f   & \wtl X \ar[l]_i \ar[dl]^{\wtl f}& X_\eta \ar[l]_j \ar[d]
 \\ S'\ar[ur]^h   \ar[r]^g & S & &\Spec\big(\ka(\eta)\big) \ar[ll]
}}\end{equation}
We assume that 
\begin{equation}\label{MB1}
  \text{$\wtl f = f \circ i \co \wtl X \to S$ is a flat morphism.} 
\end{equation}

\begin{inparaenum}[\rm (a)]%\setcounter{enumi}[1]
\item \label{MBa} For $x\in X$, the following are equivalent : 
\end{inparaenum}
\begin{enumerate}[label={\rm (\roman*)}] 
 \item\label{MBai} $x\in \wtl X$, 
 \item \label{MBaii} $x$ is a specialisation of some point in $X_\eta$, 
 \item \label{MBaiii} there exists a flat morphism $g \co S'\to S$ and an $S$--morphism $h \co S'\to X$ such that $x$ is in the image of $h$.
 \end{enumerate}
\sm  
  
\begin{inparaenum}[\rm (a)]\setcounter{enumi}{1}
\item \label{MBc} Suppose there exists a faithfully flat morphism $g\co S'\to S$ and a morphism $h \co S' \to X$ such that $g = f \circ h$. Furthermore assume $\dim X_\eta < \infty$, e.g., assume $X_\eta$ irreducible.   Then 
 \begin{equation}\label{MBc1}  \dim X_s \ge \dim \wtl X_s \ge \dim X_\eta
 \end{equation}   
holds for all $s\in S$. 
\sm 

\item \label{MBb} Suppose that $S$ is locally noetherian. We consider the following conditions \ref{MBba}, \ref{MBbb} and \ref{MBbc}: 
    \end{inparaenum}
\begin{enumerate}[label={\rm (\greek*)}] 
\item \label{MBba} $i \co \wtl X \to X$ is an isomorphism of $S$--schemes, 
    
\item \label{MBbb} for every $s\in S$, $s\ne \eta$, the $\ka(s)$--morphism $i_s \co \wtl X_s \to X_s$ is an isomorphism, 

\item \label{MBbc} $\dim \wtl X_s = \dim X_s$ holds for all $s\in S$, $s\ne \eta$. 
\end{enumerate}
Then 
\[ \text{\ref{MBba} $\iff$ \ref{MBbb} $\implies$ \ref{MBbc}.}
\] 
If all fibres $X_s$, $s\ne \eta$, are integral, then \ref{MBba}, \ref{MBbb} and \ref{MBbc} are all equivalent. 
\end{lem}

\begin{proof} \eqref{MBa} The implication \ref{MBaii} $\implies$ \ref{MBai} is trivial, 
\lv{%%%%%%%%%%   lv start
Suppose $x\in \ol{\{ y \} }$, $y\in X_\eta$; then by continuity of $j$ we have $j(x) \in j \big(\ol{\{ y \} } \big)  \subset  \ol{ \{ j(y)\}} \subset j(X_\eta)\subset \wtl X$
}%%%% lv end  
and, by \eqref{MB0},  the implication \ref{MBai} $\implies$ \ref{MBaiii} holds by taking $S'=\wtl X$ and $h=i$. Assume \ref{MBaiii}. We then know that there exists $s'\in S'$ such that $x = h(s')$. Recall \cite[Tag 03HV]{Stacks} that generizations lift along the given flat morphism $g \co S' \to S$. Since $\eta$ is a generization of $s=f(x)$, there exists a generization $t'\in S'$ of $s'\in S'$ with $g(t') = \eta$. By continuity of $h$, we get $x = h(s') \in h\big( \ol{ \{t'\} } \big) \subset \ol{ \{ h(t')\}}$, while $\eta = g(t') = f\big( h(t')\big)$ shows $h(t') \in X_\eta$. \sm 

\eqref{MBc} By base change, $X_\eta\to \Spec\big( \ka(\eta)\big)$ is locally of finite type. Hence, if $X_\eta$ is irreducible, it is finite-dimensional, \cite[IV$_2$, 4.1.1]{EGA}. Put $d = \dim X_\eta$. 

According to the upper semi-continuity of fibre dimensions \cite[IV$_3$, 13.1.3]{EGA}, applied to the structure morphism $\wtl f: \wtl X \to S$, the set
\begin{equation*}
\Sigma= \bigl\{ y \in \wtl X :  \dim\bigl(  \wtl f^{-1}( \wtl f(y)) \bigr) \geq d \bigr\}
\end{equation*} 
is closed in $X$ and is therefore stable under specializations \cite[Tag 0062]{Stacks}. 

Pick $s\in S$. Since $g \co S' \to S$ is surjective, there exists $s'\in S'$ such that $g(s') = s$.  We then have $s= g(s') = f \big(h(s')\big)$, in particular $x=h(s) \in X_s$. Thus \ref{MBaiii} of \eqref{MBa} holds. By \ref{MBai} and \ref{MBaii} we get $x\in \wtl X$ and $x\in \ol{ \{ x_\eta\}}$ for some $x_\eta \in X_\eta$. Since $x_\eta \in \wtl f{}\me\big( \wtl f (x_\eta)\big) \in \Sigma$ and $\Sigma$ is closed under specializations, we conclude $x\in \Sigma$, proving $\dim \wtl X_s \ge \dim X_\eta$. 
\sm 

\eqref{MBb} We first show that $i_\eta \co \wtl X_\eta \to X_\eta$ is an isomorphism (this holds without the assumption that $S$ is locally noetherian). Indeed, the closed immersion $i$ is quasi-compact. %\cite[10.2(2)]{GW}
Therefore, by \ref{sim}\eqref{simb}, %\cite[14.6]{GW}
the schematic image commutes with flat base change. In particular, we can apply base change to the canonical morphism $ E = \Spec\big( \ka(\eta)\big) \to S$,  which is flat because $S$ is an integral scheme. %\cite[3.29]{GW} 
Thus, we get $(\wtl X)_\eta = \Ima(j) \times_S E \cong \Ima(j_\eta)$. Since $X_\eta = (X_\eta)_\eta \to \Ima(j_\eta) \to X_\eta$, we have $\Ima(j_\eta) = X_\eta$ and therefore $i_\eta$ is an isomorphism.  

We next prove that \ref{MBba} and \ref{MBbb} are equivalent. The closed immersion $i$ is of finite type \cite[Tag 01T5]{Stacks}. Hence the structure morphism $\wtl f = f\circ i \co \wtl X \to S$ is locally of finite type \cite[Tag 01T3]{Stacks}. It follows that both $\wtl X$ and $X$ are locally of finite presentation \cite[Tag 01TX]{Stacks} as $S$--schemes. Thus, we can apply the Fibrewise Isomorphism Criterion \cite[IV$_4$, 17.9.5]{EGA}) %\ref{ag}\eqref{ag-d}
to  conclude that $i$ is an isomorphism if and only if all morphisms $i_s \co \wtl X_s \to X_s$, $s\in S$ are isomorphisms. Since we have already shown that $i_\eta$ is an isomorphism, this finishes the proof of \ref{MBba} $\iff$ \ref{MBbb}. 

Clearly, \ref{MBbb} $\implies$ \ref{MBbc}. Conversely, integrality of the fibres implies $\dim X_s< \infty$ for all $s\ne \eta$, cf.\ the proof of \eqref{MBc}. Thus \ref{MBbc} $\implies$ \ref{MBbb}  follows from  \cite[Cor.~5.8]{GW},  applied to the closed immersion $i_s$.  \end{proof}
\ms

We combine the parts of Lemma~\ref{MB} in the following result. 

\begin{prop}\label{MBco}
  Let $S$ be a locally noetherian integral scheme and let $f\co X \to S$ be a morphism locally of finite type. We use the notation of Lemma~{\rm \ref{MB}}  and assume
\begin{enumerate}[label={\rm (\roman*)}]
  \item \label{MBcoi} $\wtl f \co \wtl X \to S$ is flat, 
  
  \item\label{MBcoii} there exists a faithfully flat morphism $g \co S'\to S$ and a morphism $h\co S'\to X$ such that $f \circ h = g$, 
  % w \footnote{AP Fix}
      
 \item \label{MBcoiv} all fibres $X_s$, $s\ne \eta$, are integral, 
     
 \item\label{MBcoiii} $\dim X_s \le \dim X_\eta< \infty$ for all $s\in S$. 
  
\end{enumerate}  
Then $i$ is an isomorphism, allowing us to identify $\wtl X = X$, and $f$ is 
faithfully flat. Moreover, if all $f_s$, $s\in S$, are smooth, then $f$ is smooth. 
\end{prop}

\begin{proof} 
Because of \ref{MBcoi}, we can apply Lemma~\ref{MB}. Since $\dim X_\eta < \infty$ by \ref{MBcoiii}, assumption \ref{MBcoii} and \ref{MB}\eqref{MBc} imply that $\dim X_s\ge \dim X_s \ge \dim X_\eta$ holds for all $s\ne \eta$, hence both inequalities are equalities by \ref{MBcoiii}. By \ref{MBcoiv} we can now use \ref{MBbc} of \ref{MB}\eqref{MBb} to conclude $\wtl X \cong X$, and so $f$ is flat. That it is also surjective, follows from \ref{MBcoii}.  

Finally, if all fibres of $X$ are smooth, then smoothness of $f$ follows from the Fiberwise Smoothness Criterion \cite[IV$_4$, 17.8.2]{EGA}. \end{proof}
\sm 

Recall that $X_S(S')$ denotes the set of $S$--morphisms $S'\to X$ of $S$--schemes $S'$ and $X$.  

\begin{cor}\label{ABGpropa} Let $S$ be an locally noetherian regular irreducible scheme of dimension $\le 1$, and let $f \co X \to S$ be a morphism locally of finite type. We can and will use the notation of Lemma~{\rm \ref{MB}}, and assume
\begin{enumerate}[label={\rm (\roman*)}]
   
  \item\label{ABGpropai} there exists a faithfully flat  %flat cover 
      $S'\to S$ for which $X_S(S')\ne \emptyset$, 
  
 \item \label{ABGpropaii} all fibres $X_s$, $s\ne \eta$, are integral, 
     
 \item\label{ABGpropaiii} $\dim X_s \le \dim X_\eta< \infty$ for all $s\in S$. 
\end{enumerate}  
Then the conclusions of\/ {\rm \ref{MBco}} hold: $i$ is an isomorphism, thus $\wtl X = X$, and $f$ is faithfully flat. Moreover, if all $f_s$, $s\in S$, are smooth, then $f$ is smooth. 
\end{cor}

\begin{proof} This is a special case of Proposition~\ref{MBco}, since $\wtl f \co \wtl X \to S$ is flat by \ref{sim}\eqref{simEGA}.
\end{proof}
%% end of lv for old proposition 
%%%%%%%%%%%%%%%%%%%%%%%%%%%%%%%%%%%%%%%%%%%%%%%%%%%%%%%%%%%%%%%%%%%
\sm

\textbf{Example:} The assumptions on $S$ in \ref{ABGpropa} are fulfilled if $S=\Spec(R)$ for $R$ a DVR. In this case, we will give another flatness criterion in  \ref{ABGlemb}.  
 
\begin{cor}\label{ABGlemb} Let $R$ be a DVR. We denote  its fraction field by $K$, its residue field by $k$, and let $\ol k$ be an algebraic closure of $k$. Let  $f\co X \to \Spec(R)$ be a morphism locally of finite type. We furthermore assume that
\begin{enumerate}[label={\rm (\roman*)}]
\item \label{ABGlembi} there exists a flat $R'\in \Ralg$ for which $\ol k \in R'\mathchar45\mathbf{alg}$, 

\item \label{ABGlembii} the corresponding map $X_R(R') \to X_R(\ol k)$ is onto, see the Examples~{\rm \ref{ABGlemex}}, 
    
 \item \label{ABGlembiii} the fibre $X_k$ is integral.
\end{enumerate}
Then the conclusions of \/ {\rm \ref{MBco}} hold: $f$ is flat surjective, and $X$ is the schematic image of the canonical morphism $X_K \to X$. \end{cor}

\begin{proof} 
We first dispose of a technicality, which is folklore. 
For any $T\in \kalg$ and any $k$--scheme $Y$ we let $Y_k(T)$ be the set of $k$--morphisms $\Spec(T) \to Y$. Let $Z$ be an $R$--scheme and let $q \co Z_k \to Z$ be the canonical morphism. Then, for any $T\in \kalg$, the map 
\begin{equation}\label{ABGlemb0}
  Z_k(T) \simlgr Z_R(T), \quad \vphi \mapsto q \circ \vphi
\end{equation} 
is a bijection. Indeed, as the standard fibre product diagram below shows, any $\psi \in Z_R(T)$ factors via a unique morphism $\Spec(T) \to Z_k$: 
\[ \xymatrix@C=60pt@R=20pt{
 \Spec(T) \ar@{.>}[dr]^{\exists!}  \ar@/^0.8pc/[drr]^\psi \ar@/_0.8pc/[ddr]_\can \\
   & Z_k \ar[r]^q \ar[d] & Z \ar[d] \\
      & \Spec(k)  \ar[r]^\can  & \Spec(R) 
} \] 

As before, we abbreviate $\wtl X = \Ima( X_K \to X)$.
Since $\Spec(R') \to \Spec(R)$ is flat, \eqref{simex2} and \eqref{simEGA1} show that $\wtl X_R(R') \cong X_R(R')$. By assumption \ref{ABGlembii}, $\wtl X_R(R') \cong X_R(R') \to X_R(\ol k)$ is onto. Since $\wtl X_R(R') \to X_R(\ol k)$ factors as $\wtl X_R(R') \to \wtl X_R (\ol k) \to X_R(\ol k)$, 
\[ \xymatrix@C=50pt{
   \wtl X_R (R') \ar[r]^{\cong} \ar[d] & X_R(R') \ar[d] \\
   \wtl X_R(\ol k) \ar[r] & X_R (\ol k)
}\]
also $\wtl X(\ol k) \to X(\ol k)$ is onto. Thus, by \eqref{ABGlemb0} for $Z=\wtl X$ and $Z=X$, the map $\wtl X_k(\ol k) \to X_k(\ol k)$ is onto. 
We can now apply \cite[I, \S3, 6.10]{DG} (or \cite[Exc.~10.6]{GW}) and conclude that $\wtl X_k \to X_k$ is a surjective morphism of $k$--schemes. Since $\wtl X_k$ is a closed subscheme of $X_k$ and since $X_k$ is integral by \ref{ABGlembiii}, this implies $\wtl X_k \cong X_k$. Thus, condition \ref{MBbb} of \ref{MB}\eqref{MBb} holds, and thus $i \co \wtl X \simlgr X$ is an isomorphism. 
Hence $f$ is flat by \eqref{simEGA1}. It is also surjective because $X_k \ne \emptyset$ by integrality of $X_k$ and therefore also $X_\eta \ne \emptyset$. 
\end{proof}

\begin{leer}{\bf Examples for \ref{ABGlemb}.} \label{ABGlemex}
We use the setting of \ref{ABGlemb}:  $R$ is a DVR with fraction field $K$ and residue field $k$, and $\ol k$ is an algebraic closure of $k$.\sm

\begin{inparaenum}[(a)]
\item \label{ABGlemex-a} Suppose $R'\in \Ralg$ has  the following properties: \end{inparaenum}
\begin{enumerate}[label={\rm (\roman*)}]
  \item \label{ABGlemex-ai}  $R'$ is a local ring dominating $R$ in the sense of \ref{lem_wsur}\ref{lem_wsur-i-i},

  \item \label{ABGlemex-aii} $R'/R$ is an integral extension.
\end{enumerate}
Then $\ol k$ is an $R'$--algebra. Indeed, domination implies that  $k$ imbeds into the residue field of $R'$, which then is an algebraic extension by \ref{ABGlemex-aii} and therefore embeds into $\ol k$.
Moreover, if
\begin{enumerate}[label={\rm (\roman*)}]\setcounter{enumi}{2}

  \item \label{ABGlemex-aiii}  $R'$ is a torsion-free $R$--module,
\end{enumerate}
then $R'$ is flat as $R$--module, because over a valuation ring flatness is equivalent to torsion-freeness, see e.g.\ \cite[VI, \S3.6, Lem~1]{BAC2} or \cite[Tag 0539]{Stacks}.
Thus, whenever \ref{ABGlemex-ai}--\ref{ABGlemex-aiii} above  hold, the assumptions on $R'$ in \ref{ABGlemb} are fulfilled as soon as
\begin{equation}
  \label{ABGlemex-1} \text{$X(R') \to X(\ol k)$ is onto.}
\end{equation}
Of course, we still need $X$ to be locally of finite type and $X_k$ to be integral.%%
\sm

\begin{inparaenum}[(a)]\setcounter{enumi}{1}
\item\label{ABGlemex-h} Let $R$ be an henselian DVR and let $\ol R$ be the integral closure of $R$ in an algebraic closure $\ol K$ of $K$. Then \ref{ABGlemex-aii} and \ref{ABGlemex-aiii} of \eqref{ABGlemex-a} obviously hold for $R'=\ol R$.  But also \ref{ABGlemex-ai} is satisfied, even in a more precise form: \end{inparaenum}
\begin{enumerate}[label={\rm (\roman*)$'$}]
\item \label{ABGlemex-hi}  $\ol R$ is a henselian local domain dominating $R$; the residue field of $\ol R$ is isomorphic to $\ol k$.  \end{enumerate}
Proof of \ref{ABGlemex-hi}: We write $\ol K = \limind_{\,\la} K_\la$ as a direct limit of its subextensions $K_\la/K$ of finite degree and let $R_\la$ be the integral closure of $R$ in $K_\la$. Then $R=\limind_{\, \la} R_\la$, where each $R_\la$ is a henselian DVR, in particular a noetherian local domain, as we have seen in the proof of \ref{lem_wsur}\eqref{lem_wsur-b}. Since the transition maps are local homomorphisms, $\ol R$ is a henselian local domain by \cite[I.3, Prop.~1]{Ray-hensel}. Let $\ol \m$ be its maximal ideal.  By \eqref{ABGlemex-a}, the field $F=\ol R/\ol \m$ is an algebraic extension field of $k$. So it remains to show that $F$ is algebraically closed. To this end, let $P\in F[X]$ be a monic polynomial and let $P_0\in \ol R[X]$ be a monic lift of $P$. The roots of $P_0$ in $\ol K$ are integral over $\ol R$, hence lie in $\ol R$, since $\ol R$ is integrally closed. Thus, $P_0$ is a product of monic polynomials of degree $1$ and then the analogous fact holds for $P$, finishing the proof of \ref{ABGlemex-hi}. \lv{%%%%%%%%%%%%%%   lv starts  %%%%%%%%%%%%%%%%%%%%%%%%%%%%%%%%%%%%%%%%%%%%%%%%%
(2022-02-26) One cannot conclude that $\ol R$ is noetherian by invoking \cite[IX, App.1, Prop.~1]{BAC3} since it is not clear that the assumption $\m_\mu = R_\mu \m_\la$ for $\la \le \mu$ is fulfilled. }%%%% lv ends %%%%%%%%%%%%%%%%%%%%%%%
Thus, as in \eqref{ABGlemex-a}, only the assumptions on $X$ and \eqref{ABGlemex-1} remain.
\sm

\begin{inparaenum}[(a)]\setcounter{enumi}{2}

\item \label{ABGlemex-b} Let us describe another setting where the assumptions in \ref{ABGlemb} hold. Namely, by \cite[IX App., Cor. of Thm.~1]{BAC3}, there exists an inflation (= {\em gonflement\/} in French) $R'$ of $R$, thus in particular a local ring, whose residue field is isomorphic to $\ol k$. Hence $\ol k$ is an $R'$--ring. Since any inflation of $R$ is a faithfully flat $R$--module \cite[IX, App., Prop.~2]{BAC3}, the assumption \ref{ABGlemb} becomes \eqref{ABGlemex-1}.
\end{inparaenum}
\end{leer}

\begin{prop}[Strong surjectivity]\label{prop_strong_surj} Let $R$ be a henselian DVR and let $\ol R$ be the integral closure of $R$ in an algebraic closure of the fraction field $\Frac(R)$ of $R$. We have seen in\/ {\rm \ref{ABGlemex}\ref{ABGlemex-hi}} that $\ol R$ is a henselian local domain whose residue field $\ol k$ is an algebraic closure of the residue field of $R$.
We assume that $f\co X \to \Spec(R)$ is a morphism locally of finite type. \sm

\begin{enumerate}[label={\rm (\alph*)}]
\item \label{prop_strong_surj-a}{\rm  (Moret-Bailly)} If $f$ is universally open and surjective, e.g., faithfully flat, then $X(\ol{R}) \to X(\ol{k})$ is  surjective.

  \item \label{prop_strong_surj-b} Conversely, if $X_k$ is integral and $X(\ol{R}) \to X(\ol{k})$ is  surjective, then $f$ is faithfully flat. 
\end{enumerate}
\end{prop}

\begin{proof} \ref{prop_strong_surj-a} Our proof is variant of the proof of  \cite[IV$_4$, 14.5.8]{EGA}.

We are given  $x_0 \in X(\ol k)$. To show surjectivity of $X(\ol R) \to X(\ol k)$, it is enough to find a local ring $B$ dominating $R$ such that $B/R$ is a finite, hence integral extension and such that $x_0$ belongs to the image of $X(B) \to X(\ol{k})$. Indeed, $B$ embeds into $\ol R$ by \ref{ABGlemex}\eqref{ABGlemex-a}, making $\ol k$ a $B$--algebra, so that the map $X(B) \to X(\ol k)$ factors through $X(\ol R) \to X(\ol k)$.
%We have $B \to \ol R \to \ol k$, hence $X(B) \to X(\ol R) \to X(\ol k)$

Let $Z$ be an irreducible component of $X$ containing $x_0$. According to the implication (b) $\Rightarrow$ (e) of  \cite[IV$_3$, 14.5.6]{EGA}, %\cite[14.14]{GW}
there exists a locally closed subscheme $Z'$ of $X$ such that $x_0 \in Z'\subset Z$ and such that the restriction $f \co Z' \to \Spec(R)$ is a dominant quasi-finite morphism. The local ring  $\calO_{Z',x_0}$ is then quasi-finite over the henselian ring $R$, so that  $\calO_{Z',x_0}= B \times C$ with $B$ being finite over  $R$ and $C_k=0$ \cite[Tag 04GG(13)]{Stacks} (this commutative algebra results is based on \cite[IV$_4$, 18.2.1]{EGA}). Since $\calO_{Z',x_0}$ dominates  $R$ and $C_k = 0$, it follows that $B$ is a local ring dominating $R$. Also, by construction, $x_0$ belongs to the image of $X(B) \to X(\ol k)$.
\lv{%%%%%%%%%%%%%%%%%%%   lv starts %%%%%%%%%%%%%%%%%%%%%%%%%%%%%%%%%%%%
Let $j_{x_0} \co \Spec(\calO_{Z', x_0}) \to Z'$ be the canonical morphism \cite[(3.4.1)]{GW}. The composition
\[ \Spec(B) \subset \Spec(B) \sqcup \Spec(C) = \Spec(\calO_{Z', x_0}) \xrightarrow{j_{x_0}} Z' \subset X
\]
maps the closed point of $\Spec(B)$ (= maximal ideal of $B$) onto $x_0$.
}%%%%%%%%%%%%%%%  lv ends %%%%%%%%%%%%%%%%%%%%%%%%%%%%%%%%%%%%%%
As pointed out at the beginning of the proof, this implies that $x_0$ belongs to the image of  $X(\ol{R}) \to X(\ol k)$.

Finally, if $f$ is flat and locally of finite type (= locally of finite  presentation since $R$ is noetherian), %\cite[Tag 01TX]{Stacks}
then $f$ is universally open \cite[IV$_2$, 2.4.6]{EGA} (or \cite[Tag 01UA]{Stacks}).
%also follows from \cite[14.35]{GW} since $\Spec(R)$ is noetherian and thus finite type = finite presentation
\sm

\ref{prop_strong_surj-b} By Example \ref{ABGlemex}\eqref{ABGlemex-h}, the assumptions of Lemma~\ref{ABGlemb} are fulfilled with $R'=\ol R$. Hence, that lemma proves that $f$ is faithfully flat. 
\end{proof}

%\comments{(2022-02) I did not find Dedekind schemes in \cite{Stacks}. In \cite[\S7.13]{GW}, a Dedekind scheme is noetherian integral scheme $S$ for which $\Ga(U, \calO_S)$ is Dedekind ring for every open affine subscheme $U \subset X$. Thus, a Dedekind scheme in the sense of \cite{GW} is necessarily irreducible, while this is not the case in our definition. The definition of a Dedekind scheme in Liu's book (\cite[4.1.1]{Liu}) is more general than ours: he only assumes that $X$ is locally noetherian.}
%\pcomments{(2022-02-06) Liu also defined similarly the notion of Dedekind scheme but [GW] is more precise. There is a slightly more general notion (allowing for example disjoint sums)  used by Raynaud and al in SGA3 VI$_B$, 12.9, namely a scheme which is locally noetherian, regular, and  of dimension $\le 1$. Many statements generalize to this setting.}%% end of PG comments

\begin{leer}{\bf Dedekind schemes}\label{dede} 
We use the notion of a Dedekind scheme as defined in \cite[1.1]{BLR}: a {\em Dedekind scheme $S$ \/} is a  noetherian scheme of dimension $\le 1$, which in addition is normal, equivalently regular (recall that a noetherian ring of dimension $\le 1$ is normal if and only if it is regular). %\cite[Rem.~B.75(2)]{GW} 

This concept of a Dedekind scheme is more general than that of \cite[\S7.13]{GW}, where it is in addition assumed that $S$ be integral, equivalently irreducible. 
% integral = reduced + irreducible, but a normal scheme is reduced [St; 033K] 
However, as \eqref{dede1} below shows, the difference is not big. \sm 

\begin{inparaenum}[(a)] \item\label{dede1} 
{\em A scheme is a Dedekind scheme if and only if it is a finite disjoint union of irreducible Dedekind schemes.}  

Indeed, assume that $S$ is a Dedekind scheme. By \cite[Tag 033M]{Stacks} and its proof, the finitely many irreducible components of $S$ are normal integral schemes and $S$ is the disjoint union of them, say $S = \bigsqcup_i S_i$. Clearly, every $S_i$ is noetherian and of dimension $\le 1$. Thus, $S_i$ is an irreducible Dedekind scheme. The converse is easy to verify. 

In particular, \eqref{dede1} shows that a Dedekind scheme is irreducible if and only if it is connected if and only if it is integral. \sm  

\item\label{dede-an} 
{\em A scheme $S$ is a Dedekind scheme if and only if $S$ has a finite open covering by  spectra of Dedekind rings.} 

This follows from \eqref{dede1} and the local characterization of schemes that are integral or normal or noetherian. 
\lv{%%%%%%%%%%%%%%%%%%   lv start
Indeed, $S$ is a noetherian scheme if and only if $S$ has a finite open covering by spectra of noetherian rings \cite[I, (2.7.1)]{EGA-neu}. A scheme is normal if and only if for every affine open subscheme $\Spec(A)$ the ring $A$ is normal, if and only if it has an open covering by spectra of normal rings \cite[Tag 033J]{Stacks}. It therefore follows that a scheme is Dedekind if and only if it has a finite open covering by spectra of noetherian normal rings of dimension $\le 1$. Thus to prove the claim, we may assume that $S=\Spec(A)$ with $A$ noetherian normal of dimension $\le 1$.   Such a scheme is a finite disjoint union of normal integral schemes \cite[033M]{Stacks}, which are therefore open affine and hence of the form $\Spec(B)$ for $B$ a Dedekind ring. 
}%%%%%%%%%%%%%%   end of lv 
\end{inparaenum}
\sm

The concept of a Dedekind scheme is convenient to eliminate the irreducibility assumption of Corollary~\ref{ABGpropa}, see Proposition~\ref{AGB1}. This proposition generalizes \cite[B.1]{AGi}, 
stated in loc.\ cit.\ for $S=\Spec(R)$, $R$ a Dedekind domain, and $X$ an $R$--group scheme. In that case the assumption \ref{AGB1}\ref{AGBl-ii} is trivially fulfilled since $X$ is a group scheme and thus $X(R) \ne \emptyset$.
\end{leer}
 
\begin{prop}\label{AGB1} Let $S$ be a Dedekind scheme and let $X$ be a
$S$--scheme which
\begin{enumerate}[label={\rm (\roman*)}]

\item \label{AGBl-i} is locally of finite type, has integral fibres of the same dimension $d\ge 0$, and

\item\label{AGBl-ii} for which there exists a 
    faithfully flat morphism  $S'\to S$ such that $X_S(S') \ne \emptyset$.   
\end{enumerate}    
Then $X$ is faithfully flat. If, furthermore, all fibres are smooth, then $X$ is $S$--smooth. 
\end{prop}

\begin{proof} By \ref{dede}\eqref{dede1} we can write $S$ as a finite disjoint union of irreducible Dedekind schemes, say $S = \bigsqcup_i S_i$. If $f \co X \to $S is the structure morphism, we get $X= \bigsqcup_i X_i$ with $X_i = f\me (S_i)$, and $f_i = f|_{X_i} \co X_i \to S_i$ satisfies \ref{AGBl-i}. Similarly, let $g \co S'\to S$ be a faithfully flat morphism for which $X_S(S') \ne \emptyset$. Then $S'=\bigsqcup S'_i$ with $S'_i = g\me (S_i)$ and $g_i \co S'_i \to S_i$ satisfies \ref{AGBl-ii}. By \ref{ABGpropa}, $f_i \co X_i \to S_i$ is faithfully flat, even smooth if all fibres are so. This implies our claim. \end{proof}
\sm 

A frequently used assumption in this section was the existence of a faithfully flat morphism $S'\to S$. We finish the section with a result that weakens this assumption {\em in case $S'\to S$ is also locally of finite presentation} (which is the same as being locally of finite type in case $S$ is locally noetherian \cite[Tag 01TX]{Stacks}). We will refer to such a morphism as an ``fppf cover". More generally, we will use the term {\em fppf cover} as defined in \cite[Tag 021M]{Stacks} and {\em fpqc cover} as defined in \cite[Tag 022B or Tag 03NV]{Stacks}, which is known to yield the same covers as the ones defined in \cite[p.~28]{Vis}. 

\begin{lem}\label{lem_fpqc}
Let $S$ be a Dedekind scheme. Let $F$ be a contravariant $S$--functor locally of finite presentation in the sense of \cite[Tag 049J]{Stacks}. Then 
the following assertions are  equivalent:

\begin{enumerate}[label={\rm (\roman*)}] 

\item \label{lem_fpqc1} There exists a fpqc cover $(S_i)_{\in I}$ of $S$ such that $F(S_i) \not = \emptyset$ for each $i \in I$; 

\item \label{lem_fpqc2} There exists a fppf cover $S' \to S$  such that $F(S') \not = \emptyset$;

\item \label{lem_fpqc3} There exists an affine quasi-finite fppf cover $S'' \to S$ 
such that $F(S'') \not = \emptyset$. 
\end{enumerate}
\end{lem}
\sm 

This applies of course to the $S$-functor $h_X$ for an $S$-scheme $X$ locally of finite presentation, see for example \cite[Thm.~10.63]{GW}. We also note that the implications \ref{lem_fpqc1} $\Longleftarrow$ \ref{lem_fpqc2} $\iff$  \ref{lem_fpqc3}
are true for any scheme $S$.

\begin{proof} By \ref{dede}\eqref{dede-an} we can assume that $S=\Spec(R)$ with $R$ a Dedekind ring.   

\ref{lem_fpqc1} $\implies$ \ref{lem_fpqc2}: 
According to \cite[Tag 022E]{Stacks}, it is harmless
to assume that the $S_i$'s are affine. We write $f_i: S_i=\Spec(A_i) \to S$ and 
$A_i= \limind_{\, j \in J_i} A_{i,j}$ as its inductive limit of finitely generated $R$-subalgebras. Since $A_i$ is flat (= torsion-free) over $R$, so are the  
$A_{i,j}$'s. We put $F(T) = F\big(\Spec(T)\big)$ for $T\in \Ralg$. Then, by  assumption on $F$, we have  $F(A_i)=\limind_{j \in J_i} F(A_{i,j})$. Hence, there exists $j_i$ such that  $F(A_{i,j})$ is not empty. The structural morphism $g_i: V_i= \Spec(A_{i,j_i})\to S$ is flat and  of finite presentation (since $S$ is noetherian), so is open by \cite[Tag 01UA]{Stacks}. In particular,  $g_i(V_i)$ is open in $S$. The map $f_i:\Spec(A_i)\to S$ factors through $V_i$, so that $f_i(S_i) \subset g_i(V_i)$, whence $S= \bigcup_{i \in I}  g_i(V_i)$. Thus, the $V_i$'s form an  fppf cover of $S$. 

 Since $S$ is quasi-compact, we can select a finite subset $L \subset I$ such that 
$S= \bigcup_{i \in L} g_i(V_i)$. Then $S'=  \coprod_{i \in L} V_i$ is an affine scheme, which is a fppf cover of $S$ and satisfies $F(S') \not=\emptyset$.
\sm

\ref{lem_fpqc2} $\implies$ \ref{lem_fpqc3}:
We are given an fppf cover  $S' \to S$  such that $F(S') \not = \emptyset$.
According to Grothendieck's quasi-finite section theorem
\cite[IV$_4$, 17.16.2]{EGA}, there exist a quasi-finite affine
fppf cover $S'' \to S$ and a morphism $h\co S''\to S'$. Hence, $F$ being contravariant, there exists a map $F(S') \to F(S'')$. Since $F(S') \ne \emptyset$, this implies $F(S'') \ne \emptyset$. 
\sm 

\ref{lem_fpqc3} $\implies$ \ref{lem_fpqc1}: This is obvious since an fppf cover is also an fpqc cover \cite[Tag 022C]{Stacks}.  
\end{proof}

\end{document}